\documentclass[11pt,reqno]{amsart}
\usepackage[letterpaper,margin=1in]{geometry}
\usepackage{amsmath,amssymb,amsthm,mathtools,microtype}
\usepackage[hidelinks]{hyperref}
\hypersetup{pdftitle={Planar convexity of isotropy orbits via cohomology and common zeros},pdfauthor={Yongzhi Luan and Tin-Yau Tam}}
\newtheorem{mainthm}{Theorem}
\newtheorem{maincor}[mainthm]{Corollary}
\newtheorem{theorem}{Theorem}[section]
\newtheorem{proposition}[theorem]{Proposition}
\newtheorem{lemma}[theorem]{Lemma}
\newtheorem{corollary}[theorem]{Corollary}
\theoremstyle{remark}
\newtheorem{remark}[theorem]{Remark}
\newtheorem{example}[theorem]{Example}
\newcommand{\R}{\mathbb R}
\newcommand{\C}{\mathbb C}
\newcommand{\Z}{\mathbb Z}
\newcommand{\TT}{\mathbb T}
\newcommand{\gfr}{\mathfrak g}
\newcommand{\kfr}{\mathfrak k}
\newcommand{\pfr}{\mathfrak p}
\newcommand{\UU}{\mathrm U}
\DeclareMathOperator{\Sym}{Sym}
\DeclareMathOperator{\diag}{diag}
\DeclareMathOperator{\tr}{tr}
\DeclareMathOperator{\conv}{conv}
\title[Planar convexity of isotropy orbits]{Planar Convexity of Isotropy Orbits via Cohomology and Common Zeros}
\author{Yongzhi Luan}
\address{School of Mathematics, Shandong University, Jinan, Shandong 250100, China}
\email{luanyongzhi@email.sdu.edu.cn}
\author{Tin-Yau Tam}
\address{Department of Mathematics and Statistics, University of Nevada, Reno, NV 89557, USA}
\email{ttam@unr.edu}
\date{\today}
\begin{document}
\begin{abstract}
Let $\mathfrak{g}$ be a noncompact simple real Lie algebra with a Cartan decomposition $\mathfrak{g}=\mathfrak{k}\oplus\mathfrak{p}$. Let $K$ denote the connected Lie group corresponding to $\mathfrak{k}$. Then $K$ acts on $\mathfrak{p}$ via the adjoint representation. In this paper, we prove that for every $C\in\mathfrak{p}$ and every real linear map $L:\mathfrak{p}\to\mathbb{R}^2$, the linear image $L(K\cdot C)$ of the orbit $K\cdot C\subset\mathfrak{p}$ is convex in $\mathbb{R}^2$, except for the circle orbits when $\mathfrak{g}\cong\mathfrak{sl}_2(\mathbb{R})$ and $C\ne0$. By the two-piece property of isotropy orbits, any point omitted from the interior of the convex hull of a planar projection induces a nonzero class in the first integral cohomology of the orbit. This class vanishes directly in the non-Hermitian case, whereas in the Hermitian case we show that it vanishes by using the central circle action and a common-zero theorem of Gichev for Laplace eigenfunctions. Our approach treats classical and exceptional noncompact simple real Lie algebras uniformly and extends to real semisimple Lie algebras via Minkowski sums.
\end{abstract}
\maketitle

\tableofcontents

\section{Introduction}\label{sec:introduction}
Let $\gfr$ be a noncompact simple real Lie algebra with a Cartan
decomposition $\gfr=\kfr\oplus\pfr$. Let $G=\operatorname{Int}(\gfr)$
be the adjoint group of $\gfr$, and let $K$ denote the connected Lie
subgroup of $G$ corresponding to $\kfr$ (that is, the analytic subgroup
of $G$ with Lie algebra $\operatorname{ad}_{\gfr}(\kfr)$). Since $K$ is a maximal
compact subgroup of $G$, we know that the homogeneous quotient $G/K$ is a
Riemannian symmetric space of noncompact type (cf. \cite[Chapter~VI, Theorem~1.1, page~252]{Helgason1978}). Moreover, since
$\gfr=\kfr\oplus\pfr$, the differential of the natural projection
$\pi:G\to G/K$ at the identity element $e$ identifies $\pfr$ with the
tangent space $T_o(G/K)$ at the origin coset $o=eK$. Because
$\pi(kgk^{-1})=k\cdot\pi(g)$ for all $k\in K$ and $g\in G$,
differentiating at $g=e$ shows that under the identification
$\pfr\cong T_o(G/K)$, the adjoint representation of $K$ on $\pfr$
coincides with the isotropy representation of $K$ on $T_o(G/K)$
(cf.~\cite[Chapter~IV, Section~3, pages~209--215]{Helgason1978}).

Fix a $K$-invariant real inner product $\langle\cdot,\cdot\rangle$ on
$\pfr$, such as the restriction of the Killing form of $\gfr$. For
$C,A,B\in\pfr$, we consider the planar numerical range
\[
 W_C(A,B)=\{(\langle A,kC\rangle,\langle B,kC\rangle):k\in K\}
 \subset\R^2.
\]
Equivalently, we study the linear image $L(K\cdot C)\subset\R^2$ of the
orbit $K\cdot C\subset\pfr$ under an arbitrary real linear map
$L:\pfr\to\R^2$. Our main result is the following theorem.

\begin{mainthm}\label{conv-main}
Let $\gfr$ be a noncompact simple real Lie algebra with a Cartan
decomposition $\gfr=\kfr\oplus\pfr$, and let $K$ denote the connected
compact Lie group corresponding to $\kfr$, acting on $\pfr$ via the
adjoint representation. If $\gfr\not\cong\mathfrak{sl}_2(\R)$, then
$L(K\cdot C)$ is convex for every $C\in\pfr$ and every real linear map
$L:\pfr\to\R^2$. For $\gfr\cong\mathfrak{sl}_2(\R)$ and $C\ne0$,
the orbit $K\cdot C$ is a circle in $\pfr$; its image $L(K\cdot C)$ is
an ellipse when $\operatorname{rank}L=2$, and a line segment or
a point when $\operatorname{rank}L\le1$.
\end{mainthm}

For complex matrices, Westwick~\cite{Westwick1975} established the
convexity of the $C$-numerical range using Morse theory. Poon~\cite{Poon1980}
later gave an elementary proof based on majorization and inclusion
relations, without using Morse theory. Westwick's result was an important
motivation for the Lie-theoretic study of numerical ranges by Li and
Tam~\cite{LiTam2000}. Li and Tam proved Theorem~\ref{conv-main} for
classical simple real Lie algebras and for complex simple Lie algebras
viewed as real Lie algebras. However, the twelve noncomplex noncompact
exceptional real forms were left open in
\cite[Remark~12.1, page~169]{LiTam2000}. Here we give a unified proof
for all noncompact simple real Lie algebras without case-by-case
analysis of the individual isotropy representations.

Isotropy representations are polar, and their orbit convex hulls have
also been studied from the viewpoint of convex algebraic geometry.
Kobert and Scheiderer~\cite{KobertScheiderer2022} proved that the convex
hull of every orbit of a polar representation is a spectrahedron. Our
result is complementary: for the isotropy orbits in
Theorem~\ref{conv-main}, apart from the stated
$\mathfrak{sl}_2(\R)$ exception, a planar linear image of an orbit
already equals the corresponding image of its convex hull.

By decomposing a real semisimple Lie algebra into simple ideals, we also
get the following extension of Theorem~\ref{conv-main}. Note that
compact simple ideals have zero $\pfr$-component in the Cartan
decomposition and therefore do not affect the orbit $K\cdot C$.

\begin{maincor}\label{conv-semisimple}
Let $\gfr$ be a real semisimple Lie algebra with a Cartan decomposition
$\gfr=\kfr\oplus\pfr$, and let $K$ denote the connected compact Lie
group corresponding to $\kfr$, acting on $\pfr$ via the adjoint
representation. If $\gfr$ has no simple ideal isomorphic to
$\mathfrak{sl}_2(\R)$, then
\[
 L(K\cdot C)\subset\R^2
 \quad\text{is convex}
\]
for every $C\in\pfr$ and every real linear map $L:\pfr\to\R^2$.
More generally, write $\gfr=\bigoplus_j\gfr_j$ as a direct sum of simple
ideals, with $\pfr=\bigoplus_j\pfr_j$, and write $C=\sum_j C_j$ with
$C_j\in\pfr_j$. Then $L(K\cdot C)$ is convex whenever
\[
 \operatorname{rank}(L|_{\pfr_j})\le1
\]
for every simple ideal $\gfr_j\cong\mathfrak{sl}_2(\R)$ with $C_j\ne0$.
Conversely, if $\gfr$ has a simple ideal isomorphic to
$\mathfrak{sl}_2(\R)$, then there exist $C\in\pfr$ and a real linear map
$L:\pfr\to\R^2$ such that $L(K\cdot C)$ is not convex.
\end{maincor}

Corollary~\ref{conv-semisimple} shows that $L(K\cdot C)$ is convex for
all $C\in\pfr$ and all real linear maps $L:\pfr\to\R^2$ if and only if
$\gfr$ has no simple ideal isomorphic to $\mathfrak{sl}_2(\R)$. For a
fixed pair $(C,L)$, the rank condition
$\operatorname{rank}(L|_{\pfr_j})\le1$ is sufficient for convexity,
though not necessary: when several $\mathfrak{sl}_2(\R)$ ideals are
present, the Minkowski sum of their nonconvex elliptical images can still
fill the interior and yield a convex set (see
Corollary~\ref{conv-exceptional-filling}).

Throughout this paper, for any subset $A$ of a real vector space, we
use $\conv A$ to denote the convex hull of $A$. Our proof of
Theorem~\ref{conv-main} combines the topological two-piece property of
isotropy orbits due to Bott and Samelson~\cite{BottSamelson1958} with a
common-zero theorem for Laplace eigenfunctions due to
Gichev~\cite{Gichev2004}. Specifically, the paper is organized as
follows. In Section~\ref{sec:angular}, for a compact connected smooth
manifold $M$, we study continuous maps $F:M\to\R^2$ whose open
half-plane preimages are connected or empty, and we show that any point
omitted from the interior of the convex hull $\conv F(M)$ defines a
nonzero angular cohomology class in $H^1(M;\Z)$. In
Section~\ref{sec:zero}, we recall Gichev's common-zero theorem for
Laplace eigenfunctions and use the central circle action to show that
this cohomology obstruction vanishes for nontrivial irreducible complex
representations of compact connected semisimple Lie groups. In
Section~\ref{sec:isotropy}, we prove Theorem~\ref{conv-main} and
Corollary~\ref{conv-semisimple}, and we give applications to complex
symmetric matrices and constrained linear optimization on isotropy
orbits. Finally, in Section~\ref{sec:extensions}, we extend our
convexity theorem to general irreducible orthogonal representations with
the two-piece property, analyze Minkowski sums of circle orbits, and
explain why Theorem~\ref{conv-main} does not extend to linear maps into
$\R^n$ for $n\ge3$.

\section{Half-plane sections and angular maps}\label{sec:angular}

Let $M$ be a compact connected smooth manifold and let $F:M\to\R^2$ be a
continuous map. Assume that
\begin{equation}\label{conv-halfplane}
 F^{-1}(H)\text{ is connected or empty for every open affine
 half-plane }H\subset\R^2.
\end{equation}
Recall that an embedded submanifold $M\subset V$ of a real vector space
$V$ is said to have the \emph{two-piece property} if its intersection
with every open affine half-space of $V$ is connected or empty. If
$M\subset V$ has the two-piece property, then the restriction to $M$ of
every real linear map $L:V\to\R^2$ satisfies \eqref{conv-halfplane},
because the inverse image of an open half-plane in $\R^2$ is either
empty, $M$, or the intersection of $M$ with an open affine half-space of
$V$.

For any subset $A\subset\R^2$, we denote by $\conv A$ the convex hull of
$A$ in $\R^2$.

\begin{lemma}\label{conv-boundary}
Put $D=\conv F(M)$. Under the assumption \eqref{conv-halfplane},
$\partial D\subseteq F(M)$ if $D$ has nonempty interior in $\R^2$.
If $\dim D\le1$, then $F(M)=D$.
\end{lemma}
\begin{proof}
Since $M$ is compact and connected and $F:M\to\R^2$ is continuous, the
image $F(M)$ is a nonempty compact connected subset of $\R^2$. If
$\dim D=0$, then $D$ is a single point, and $\emptyset\ne F(M)\subseteq D$
gives $F(M)=D$. If $\dim D=1$, then $D$ is contained in an affine line
$\ell\subset\R^2$. Identifying $\ell$ with $\R$ by an affine
isomorphism, we see that $F(M)$ is a nonempty compact connected subset
of $\R$. Since every nonempty compact connected subset of $\R$ is a
closed bounded interval $[a,b]$ (by the extreme-value theorem and the
intermediate-value theorem), $F(M)$ is a closed line segment in $\ell$.
In particular, $F(M)$ is convex, and hence $D=\conv F(M)=F(M)$.

Suppose now that $D$ has nonempty interior in $\R^2$. Let
$\lambda:\R^2\to\R$ be a nonzero linear functional and set
$c=\max_M\lambda\circ F$. The sets
\[
 E_j=\overline{\{x\in M:\lambda(F(x))>c-1/j\}}
 \qquad(j\ge1)
\]
are nonempty, compact, and connected by \eqref{conv-halfplane}, and they
decrease with $j$. Hence their intersection
\[
 E=\{x\in M:\lambda(F(x))=c\}
\]
is also compact and connected. The corresponding exposed face of $D$ is
$D\cap\{\lambda=c\}=\conv F(E)$. Since $F(E)$ is a compact connected
subset of the line $\{\lambda=c\}$, the one-dimensional case proved
above gives $\conv F(E)=F(E)\subseteq F(M)$. Since every boundary point
of a planar compact convex set with nonempty interior lies on a
supporting line, it follows that $\partial D\subseteq F(M)$.
\end{proof}

For $p\notin F(M)$, we define the angular map
\[
 \theta_p:M\longrightarrow S^1,
 \qquad \theta_p(x)=\frac{F(x)-p}{|F(x)-p|}.
\]
Under the standard identification $[M,S^1]\cong H^1(M;\Z)$, we denote by
$[\theta_p]\in H^1(M;\Z)$ the homotopy class of $\theta_p$. Recall that
$[\theta_p]=0$ if and only if there exists a continuous function
$\alpha:M\to\R$ such that $\theta_p=e^{i\alpha}$.

\begin{lemma}\label{conv-angular}
Under the assumption \eqref{conv-halfplane}, if
$p\in\operatorname{int}D\setminus F(M)$, then $[\theta_p]\ne0$ in
$H^1(M;\Z)$.
\end{lemma}
\begin{proof}
Suppose toward a contradiction that $[\theta_p]=0$, so that
$\theta_p=e^{i\alpha}$ for a continuous function $\alpha:M\to\R$. Since
$M$ is compact and connected, we have $\alpha(M)=[a,b]$ for some
$a\le b$. If $b-a>\pi$, we can choose $s,t\in[a,b]$ with
$\pi<t-s<2\pi$. Let $U\subset S^1$ be the open semicircle centered at
$e^{i((s+t)/2+\pi)}$, which contains $e^{is}$ and $e^{it}$ but not
$e^{i(s+t)/2}$. Then $\theta_p^{-1}(U)$ is the preimage under $F$ of an
open affine half-plane whose boundary line passes through $p$. By
\eqref{conv-halfplane}, $\theta_p^{-1}(U)$ is connected, so its image
$\alpha(\theta_p^{-1}(U))$ is an interval containing $s$ and $t$. By the
intermediate-value theorem, $(s+t)/2\in\alpha(\theta_p^{-1}(U))$, which
contradicts $e^{i(s+t)/2}\notin U$.

Therefore, we must have $b-a\le\pi$. This implies that $\theta_p(M)$ is
contained in a closed semicircle, and hence $F(M)$ is contained in a
closed half-plane whose boundary line passes through $p$. This
contradicts $p\in\operatorname{int}D$.
\end{proof}

\begin{proposition}\label{conv-topological}
Under the assumption \eqref{conv-halfplane}, $F(M)$ is convex if and
only if $[\theta_p]=0$ for every $p\notin F(M)$. In particular, if
$H^1(M;\Z)=0$, then $F(M)$ is convex.
\end{proposition}
\begin{proof}
If $F(M)$ is convex and $p\notin F(M)$, we can strictly separate $p$
from the compact convex set $F(M)$ by a line. Then $\theta_p(M)$ lies in
an open semicircle, which is contractible, so $[\theta_p]=0$.
Conversely, if $F(M)$ is not convex, then by Lemma~\ref{conv-boundary}
there exists a point $p\in\operatorname{int}D\setminus F(M)$, and
Lemma~\ref{conv-angular} shows that $[\theta_p]\ne0$.
\end{proof}

\section{The common-zero method and the scalar circle}\label{sec:zero}

To treat the Hermitian case of Theorem~\ref{conv-main}, we use a
common-zero theorem for Laplace eigenfunctions due to
Gichev~\cite[Theorem~1 and Corollary~2, pages~202--204]{Gichev2004}. For
completeness, we include a short proof below.

\begin{lemma}\label{conv-eigenzero}
Let $N$ be a compact connected Riemannian manifold without boundary
such that $H^1(N;\R)=0$. If a smooth function $f:N\to\C$ satisfies
$\Delta f=-\lambda f$ with $\lambda>0$, then $f$ has a zero.
\end{lemma}
\begin{proof}
Suppose that $f$ has no zero on $N$. Write $f=u+iv$, where
$u=\operatorname{Re}f$ and $v=\operatorname{Im}f$ are smooth real-valued
functions on $N$, and define the modulus function $r:N\to\R$ by
\[
 r(x)=|f(x)|=\sqrt{u(x)^2+v(x)^2}\qquad(x\in N).
\]
Since $f(x)\ne0$ for all $x\in N$, we have $u(x)^2+v(x)^2>0$ everywhere
on $N$. Because the square-root function is smooth on $(0,\infty)$, the
function $r=|f|$ is smooth and strictly positive on $N$. Moreover, the
real one-form
\[
 \operatorname{Im}\Bigl(\frac{df}{f}\Bigr)=\frac{u\,dv-v\,du}{u^2+v^2}
\]
is smooth and closed on $N$ (locally, it is the differential of a smooth
branch of $\arg f$). Since $H^1(N;\R)=0$, this one-form is exact, so
there exists a smooth real-valued function $\alpha:N\to\R$ such that
$f=r e^{i\alpha}$. Taking the imaginary part of the equation
$\Delta f=-\lambda f$, we have
\[
 2\langle\nabla r,\nabla\alpha\rangle+r\Delta\alpha=0,
 \qquad\text{and hence}\qquad
 \operatorname{div}(r^2\nabla\alpha)=0.
\]
Multiplying by $\alpha$ and integrating by parts over $N$, we get
\[
 0=-\int_N r^2|\nabla\alpha|^2\,d\mathrm{vol}.
\]
Since $r>0$, it follows that $\alpha$ is constant. The real part of the
eigenvalue equation then reduces to $\Delta r=-\lambda r$. Integrating
both sides over $N$ gives
\[
 0=\int_N \Delta r\,d\mathrm{vol}=-\lambda\int_N r\,d\mathrm{vol}<0,
\]
which is a contradiction.
\end{proof}

\begin{lemma}[{\cite[Corollary~2, page~204]{Gichev2004}}]\label{conv-repzero}
Let $S$ be a compact connected semisimple Lie group, and let $E$ be
a nontrivial irreducible finite-dimensional complex representation
of $S$. For every $C\in E$ and every real linear map $L:E\to\C$, there
exists $s\in S$ such that $L(sC)=0$.
\end{lemma}
\begin{proof}
Fix an $S$-invariant Hermitian inner product on $E$ and a bi-invariant
Riemannian metric on $S$. Let $X_1,\ldots,X_d$ be an orthonormal basis
of the Lie algebra of $S$, and let $d\rho$ denote the derived
representation of the Lie algebra on $E$. By Schur's lemma, the Casimir
operator
\[
 -\sum_{j=1}^d d\rho(X_j)^2
\]
acts on $E$ as a scalar $cI_E$. By unitarity, we have $c\ge0$. If $c=0$,
then $\sum_j\|d\rho(X_j)v\|^2=0$ for all $v\in E$, which implies that
$d\rho=0$ and hence the representation of $S$ is trivial because $S$ is
connected, a contradiction. Thus $c>0$.

For a fixed $C\in E$, define $f:S\to\C$ by $f(s)=L(sC)$. Since the
Laplacian on $S$ is given by $\Delta=\sum_{j=1}^d X_j^2$ and $L$ is
real linear with $c\in\R$, we get $\Delta f=-cf$. Because $S$ is a
compact connected semisimple Lie group, its fundamental group is finite
and thus $H^1(S;\R)=0$. By Lemma~\ref{conv-eigenzero}, there exists
$s\in S$ such that $f(s)=L(sC)=0$.
\end{proof}

\begin{proposition}\label{conv-circle}
Let $S$ and $E$ be as in Lemma~\ref{conv-repzero}, and let
$\TT$ act on $E$ by complex scalars. If $M=\TT S\cdot C$ has
the two-piece property, then $L(M)$ is convex for every real linear map
$L:E\to\R^2$.
\end{proposition}
\begin{proof}
If $C=0$, the result is clear. Assume $C\ne0$ and identify $\R^2$ with
$\C$. By Lemma~\ref{conv-repzero}, there exists $v=s_0C\in S\cdot C$
such that $L(v)=0$. Since $L$ is real linear, for all $t\in\R$ we have
\[
 L(e^{it}v)=\cos t\,L(v)+\sin t\,L(iv)
           =\sin t\,L(iv).
\]
Hence the circle $\TT\cdot v$ is mapped by $L$ onto a line segment.

Let $P=S\times\TT$ and consider the orbit map $q:P\to M$ defined by
$q(s,z)=zsC$. Since the stabilizer $P_C$ is a compact Lie group, it has
finitely many connected components. By the homotopy exact sequence of
the fibration $P_C\to P\to M$, the subgroup $q_*(\pi_1(P))$ has finite
index in $\pi_1(M)$. Since $\Z$ is torsion-free, any homomorphism from
$\pi_1(M)$ to $\Z$ that vanishes on $q_*(\pi_1(P))$ must vanish on all
of $\pi_1(M)$. Consequently, the induced homomorphism
\[
 q^*:H^1(M;\Z)\longrightarrow H^1(P;\Z)
\]
is injective. Moreover, since $H^1(S;\Z)=0$, we have
$H^1(P;\Z)\cong H^1(\TT;\Z)\cong\Z$, and any class in $H^1(P;\Z)$ is
determined by its restriction to the slice $\{s_0\}\times\TT$.

Now let $p\notin L(M)$. The loop $t\mapsto L(e^{it}v)$ moves along a
line segment that does not contain $p$, so its winding number around $p$
is zero. This means that the restriction of $q^*[\theta_p]$ to
$\{s_0\}\times\TT$ is zero. Therefore $q^*[\theta_p]=0$, and the
injectivity of $q^*$ gives $[\theta_p]=0$. By
Proposition~\ref{conv-topological}, we conclude that $L(M)$ is convex.
\end{proof}

\section{Planar convexity of isotropy orbits}\label{sec:isotropy}

In this section, we prove Theorem~\ref{conv-main} and
Corollary~\ref{conv-semisimple}. By a classical theorem of Bott and
Samelson~\cite[Theorem~I$'$, page~965]{BottSamelson1958} (see also
Gorodski~\cite[Sections~2--3, pages~62--65]{Gorodski2006}), every orbit
$M=K\cdot C\subset\pfr$ of the isotropy representation of a Riemannian
symmetric space is taut in $\pfr$. In particular, for any open affine
half-space $H\subset\pfr$, the intersection $M\cap H$ is an increasing
union of closed metric balls (spherical caps) in $M$ whose zeroth
\v{C}ech homology injects into $H_0(M;\Z_2)\cong\Z_2$. Thus every
isotropy orbit $M=K\cdot C$ (including singular orbits) has the
two-piece property, and the restriction of any real linear map
$L:\pfr\to\R^2$ to $M$ satisfies \eqref{conv-halfplane}.

\begin{proof}[Proof of Theorem~\ref{conv-main}]
Let $M=K\cdot C$. As noted above, $M$ has the two-piece property.
We divide the proof into two cases according to the center of $\kfr$.

\emph{Case 1: $\mathfrak{z}(\kfr)=0$ (the non-Hermitian case).}
In this case, $\kfr$ is semisimple, so the connected compact Lie group
$K$ is semisimple. Write $M=K/K_C$, where $K_C$ is the stabilizer of $C$
in $K$. From the homotopy exact sequence
\[
 \pi_1(K)\longrightarrow\pi_1(M)
 \longrightarrow\pi_0(K_C)\longrightarrow1,
\]
we see that $\pi_1(M)$ is finite, because $\pi_1(K)$ is finite for the
compact semisimple group $K$ and $\pi_0(K_C)$ is finite for the compact
Lie group $K_C$. Therefore, $H^1(M;\Z)=\operatorname{Hom}(\pi_1(M),\Z)=0$.
By Proposition~\ref{conv-topological}, $L(M)$ is convex.

\emph{Case 2: $\mathfrak{z}(\kfr)\ne0$ (the Hermitian case).}
By the structure theory of irreducible Hermitian symmetric spaces (see
\cite[Chapter~VIII, Section~6, pages~514--518]{Helgason1978} and
\cite[Chapter~VII, Section~9]{Knapp2002}), $\dim\mathfrak{z}(\kfr)=1$,
and there is a $K$-invariant complex structure $J$ on $\pfr$ such that
the connected center $Z(K)^0\cong\TT$ acts on $\pfr$ by unit complex
scalars. Let $S=[K,K]$ be the derived subgroup of $K$. Then $S$ is a
compact connected semisimple Lie group, and $K=\TT S$.

Since $\gfr$ is simple, the adjoint action of $K$ on $\pfr$ is
irreducible over $\R$. Hence $E=(\pfr,J)$ is irreducible as a complex
representation of $K$. Moreover, $E$ remains irreducible when restricted
to $S$, because any complex $S$-invariant subspace is also invariant
under the scalar action of $\TT$ and hence under $K=\TT S$. If
$\gfr\not\cong\mathfrak{sl}_2(\R)$, then $\dim_\C E>1$, so $E$ is a
nontrivial irreducible complex representation of $S$. By
Proposition~\ref{conv-circle}, $L(M)=L(\TT S\cdot C)$ is convex.

Finally, if $\gfr\cong\mathfrak{sl}_2(\R)$, then $\dim_\R\pfr=2$ and
$K\cong\mathrm{SO}(2)$ acts on $\pfr\cong\R^2$ by rotations. Thus, for
$C\ne0$, the orbit $K\cdot C$ is a circle in $\pfr$, and its image under
a real linear map $L:\pfr\to\R^2$ is an ellipse when
$\operatorname{rank}L=2$, and a line segment or a point when
$\operatorname{rank}L\le1$.
\end{proof}

\begin{proof}[Proof of Corollary~\ref{conv-semisimple}]
Let $\gfr=\bigoplus_j\gfr_j$ be the decomposition of $\gfr$ into simple
ideals, with corresponding Cartan decompositions
$\gfr_j=\kfr_j\oplus\pfr_j$. Any compact simple ideal has $\pfr_j=0$ and
contributes zero to $\pfr$. For each noncompact simple ideal $\gfr_j$,
let $K_j$ denote the connected compact Lie group corresponding to
$\kfr_j$, acting on $\pfr_j$ via the adjoint representation. Since the
simple ideals commute, we have
\[
 K\cdot C=\prod_j(K_j\cdot C_j),\qquad
 L(K\cdot C)=\sum_j L_j(K_j\cdot C_j),
 \quad L_j=L|_{\pfr_j},
\]
where the sum on the right-hand side is the Minkowski sum in $\R^2$. If
$\gfr_j\not\cong\mathfrak{sl}_2(\R)$, then $L_j(K_j\cdot C_j)$ is convex
by Theorem~\ref{conv-main}. If $\gfr_j\cong\mathfrak{sl}_2(\R)$, then by
Theorem~\ref{conv-main}, $K_j\cdot C_j$ is a point when $C_j=0$, and its
image $L_j(K_j\cdot C_j)$ is a line segment or a point when $C_j\ne0$
and $\operatorname{rank}L_j\le1$. Since the Minkowski sum of convex sets
in $\R^2$ is convex, $L(K\cdot C)$ is convex.

Conversely, suppose that $\gfr$ has a simple ideal $\gfr_{j_0}$
isomorphic to $\mathfrak{sl}_2(\R)$. Choose $C\in\pfr_{j_0}$ with
$C\ne0$ (having zero components in all other simple ideals), and let
$L:\pfr\to\R^2$ be a linear map that restricts to a linear isomorphism
from $\pfr_{j_0}$ onto $\R^2$ and vanishes on $\pfr_j$ for all
$j\ne j_0$. Then $L(K\cdot C)=L_{j_0}(K_{j_0}\cdot C)$ is an ellipse in $\R^2$, which is not convex.
\end{proof}

\begin{remark}\label{conv-sl2-explanation}
We can illustrate the exception $\gfr=\mathfrak{sl}_2(\R)$ explicitly.
In the standard Cartan decomposition $\mathfrak{sl}_2(\R)=\kfr\oplus\pfr$,
we have
\[
 \kfr=\mathfrak{so}(2),\qquad
 \pfr=\left\{
 \begin{pmatrix}x&y\\y&-x\end{pmatrix}:x,y\in\R
 \right\}.
\]
For
$k_\theta=\left(\begin{smallmatrix}\cos\theta&-\sin\theta\\
\sin\theta&\cos\theta\end{smallmatrix}\right)\in\mathrm{SO}(2)$,
the adjoint action $\operatorname{Ad}(k_\theta)$ rotates the coordinates
$(x,y)$ by angle $2\theta$. Taking $C=\operatorname{diag}(1,-1)\in\pfr$,
we get the orbit
\[
 K\cdot C=
 \left\{\begin{pmatrix}
 \cos2\theta&\sin2\theta\\
 \sin2\theta&-\cos2\theta
 \end{pmatrix}:\theta\in\R\right\},
\]
whose image under the linear isomorphism
$L\left(\begin{smallmatrix}x&y\\y&-x\end{smallmatrix}\right)=(x,y)$ is
the unit circle $S^1\subset\R^2$. Here $K\cdot C\cong S^1$ satisfies the
two-piece property, but $H^1(K\cdot C;\Z)\cong\Z$ is nonzero and the
derived subgroup $S=[K,K]$ is trivial, so the hypotheses of
Proposition~\ref{conv-topological} and Proposition~\ref{conv-circle} are
not satisfied.

This exceptional two-dimensional orbit is consistent with Brickman's
theorem~\cite{Brickman1961}. For $n\ge3$ and real symmetric $n\times n$
matrices $A,B$, Brickman proved that
\[
 \{(x^{\top}Ax,x^{\top}Bx):x\in\R^n,\ x^{\top}x=1\}
\]
is convex. In dimension two, by contrast, a pair of traceless symmetric
matrices can give a circle, as the example above shows.
\end{remark}

As an immediate application of Theorem~\ref{conv-main}, we consider the
unitary congruence orbit of a complex symmetric matrix.

\begin{corollary}\label{conv-symmetric}
Let $n\ge2$ and let $S=\operatorname{diag}(s_1,\ldots,s_n)$ with
$s_j\ge0$. For every real linear map $L:\Sym_n(\C)\to\R^2$, the set
\[
 \{L(U^{\top}SU):U\in\UU(n)\}
\]
is convex.
\end{corollary}
\begin{proof}
The unitary congruence action of $\UU(n)$ on $\Sym_n(\C)$ is
equivalent to the isotropy representation of the noncompact simple real
Lie algebra $\mathfrak{sp}(n,\R)$ (note that for $n\ge2$,
$\mathfrak{sp}(n,\R)\not\cong\mathfrak{sl}_2(\R)$). Thus the result
follows from Theorem~\ref{conv-main}.
\end{proof}

Note that restricting to two real coordinates in
Corollary~\ref{conv-symmetric} is essential. Even for $n=2$, the set of
all diagonals of matrices in the congruence orbit of $I_2$ is not
convex: it contains $(1,1)$ and $(1,-1)$, whereas their midpoint $(1,0)$
cannot be the diagonal of a symmetric unitary matrix (since any unitary
matrix with first diagonal entry $1$ must have zero off-diagonal entries
in the first row and column, so its second diagonal entry must have
modulus $1$).

Theorem~\ref{conv-main} also implies that linear optimization problems
on an isotropy orbit with a single linear constraint have zero
convex-relaxation gap.

\begin{corollary}\label{conv-optimization}
Let $M=K\cdot C\subset\pfr$ be an isotropy orbit as in
Theorem~\ref{conv-main} with $\gfr\not\cong\mathfrak{sl}_2(\R)$.
Let $a,b:\pfr\to\R$ be real linear functionals, and let $I\subset\R$ be
a closed interval, possibly unbounded. If the feasible set
$\{x\in\conv M:b(x)\in I\}$ is nonempty, then
\[
 \max_{x\in M,\ b(x)\in I}a(x)
 =\max_{y\in\conv M,\ b(y)\in I}a(y).
\]
Moreover, for every $y\in\conv M$ with $b(y)\in I$, there exists
$x\in M$ such that $a(x)=a(y)$ and $b(x)=b(y)$.
\end{corollary}
\begin{proof}
Define $L:\pfr\to\R^2$ by $L(x)=(a(x),b(x))$. By
Theorem~\ref{conv-main}, $L(M)$ is convex, so we have
$L(\conv M)=\conv L(M)=L(M)$. This gives the second assertion, and the
equality of the maxima follows from the compactness of $M$.
\end{proof}

\section{Extensions, examples, and related results}\label{sec:extensions}

The proof of Theorem~\ref{conv-main} relies only on the decomposition of
$K$ into its connected center and its derived semisimple subgroup. The
same argument therefore extends to general irreducible orthogonal
representations of compact connected Lie groups.

Before stating this extension, we recall the definition of a taut
submanifold and a taut representation (see Carter and
West~\cite[page~707]{CarterWest1972}, Gorodski and
Thorbergsson~\cite[page~187]{GorodskiThorbergsson2003}, and
Gorodski~\cite[Sections~2--3, pages~62--65]{Gorodski2006}). Let $V$ be a
finite-dimensional real Euclidean space. A compact connected submanifold
$M\subset V$ is called \emph{taut} if for every point $p\in V$ and every
radius $r\ge0$, the inclusion of the closed Euclidean ball section
$B_M(p,r)=\{x\in M:\|x-p\|\le r\}$ into $M$ induces an injective
homomorphism
\[
 H_*(B_M(p,r);\Z_2)\longrightarrow H_*(M;\Z_2)
\]
on \v{C}ech homology with $\Z_2$-coefficients (equivalently, every Morse
distance function $L_p(x)=\|x-p\|^2$ on $M$ has the minimum number of
critical points required by the $\Z_2$-Betti numbers of $M$). An
orthogonal representation of a compact Lie group $K$ on $V$ is called a
\emph{taut representation} if every orbit $K\cdot C\subset V$ ($C\in V$)
is a taut submanifold of $V$. Since every orbit $M=K\cdot C$ of an
orthogonal representation lies on a sphere centered at the origin in
$V$, its intersection with any closed affine half-space of $V$ coincides
with a closed ball section $B_M(p,r)$, and its intersection with any
open affine half-space is an increasing union of such closed ball
sections. In particular, the injectivity on $H_0(\cdot\,;\Z_2)$ implies
that every orbit of a taut representation has the two-piece property.

\begin{theorem}\label{conv-irreducible-general}
Let a compact connected Lie group $K$ act orthogonally and
irreducibly on a real vector space $V$ with $\dim_\R V>2$.
If the orbit $M=K\cdot C$ has the two-piece property, then
$L(M)$ is convex for every real linear map $L:V\to\R^2$.
In particular, this conclusion holds for every orbit of an
irreducible taut representation in real dimension greater than two.
\end{theorem}
\begin{proof}
Replacing $K$ by its image in $\mathrm O(V)$, we may assume that the
action is effective. Put $S=[K,K]$ and $T=Z(K)^0$; then $K=TS$, where
$S$ is compact, connected, and semisimple. If $T$ is trivial, the same
homotopy exact sequence argument as in Case~1 of the proof of
Theorem~\ref{conv-main} shows that $H^1(M;\Z)=0$, and the result follows
from Proposition~\ref{conv-topological}.

Suppose now that $T$ is nontrivial. Since $T$ is central in $K$, the
real $T$-isotypic components of $V$ are $K$-invariant. By the
irreducibility of $V$ over $\R$, there is only one real $T$-isotypic
component. Since the action of $K$ is effective and $T$ is a nontrivial
connected torus, $T$ acts nontrivially on $V$. Consequently, on the
complexification of $V$, the weights of $T$ consist of a single
conjugate pair $\chi,\overline\chi$ for a nontrivial character
$\chi:T\to\TT$ with $\chi(T)=\TT$. This defines a $K$-invariant complex
structure $J$ on $V$ such that $T$ acts on $E=(V,J)$ by the complex
scalars $\chi(t)$.

Since $V$ is irreducible as a real $K$-module, $E=(V,J)$ is irreducible
as a complex $K$-module, and hence also irreducible as a complex
$S$-module (since any complex $S$-invariant subspace is invariant under
the scalar action of $T$). Because $\dim_\R V>2$, we have
$\dim_\C E>1$, so the irreducible $S$-representation on $E$ is
nontrivial. Thus $M=\TT S\cdot C$, and the convexity of $L(M)$ follows
from Proposition~\ref{conv-circle}.
\end{proof}

\begin{remark}
Taut irreducible representations of compact Lie groups were classified
up to orbit equivalence by Gorodski and
Thorbergsson~\cite{GorodskiThorbergsson2003}.
Theorem~\ref{conv-irreducible-general} provides a direct proof of the
planar convexity of their orbits (and of any individual orbit with the
two-piece property in an irreducible orthogonal representation of
dimension greater than two) without using the classification. Its proof
combines the common-zero theorem of Gichev, standard representation
theory, and the two-piece property.
\end{remark}

We note that the two-piece property alone does not imply planar
convexity when the representation is reducible. For instance, consider
the product torus $M=S^1\times S^1\subset\C^2$, which is an isotropy
orbit of the reducible symmetric space associated with
$\mathfrak{sl}_2(\R)\oplus\mathfrak{sl}_2(\R)$ (see
Example~\ref{conv-so22}). Any real linear functional on $\C^2$ restricts
to $M$ as a sum of two independent cosine functions, so its strict
superlevel sets on $M$ are connected and $M$ has the two-piece property.
However, for $0<r<1$, the real linear map $L(z,w)=z+rw$ maps $M$ onto
the annulus
\[
 \{u\in\C:1-r\le|u|\le1+r\},
\]
which is not convex.

\subsection{Minkowski sums of projected circle orbits}\label{subsec:exceptional-sums}

For subsets $A,B\subset\R^2$, their Minkowski sum is defined by
$A+B=\{a+b:a\in A,\ b\in B\}$. Sums of more than two sets are
defined in the same way.

In the semisimple setting of Corollary~\ref{conv-semisimple}, the
individual $\mathfrak{sl}_2(\R)$ orbits are circles, but their rank-two
linear images are ellipses. Even when several factors have such
nonconvex elliptical images, their Minkowski sum can be convex if the
ellipses together fill the interior. We first study circles and then
apply the result to ellipses that are images of circles under a common
invertible real linear map. We identify $\R^2$ with $\C$ and denote the
unit circle by $S^1$.

\begin{lemma}\label{conv-circle-sums}
Let $r_1,\ldots,r_q>0$ with $q\ge1$, and set $R=\sum_{j=1}^q r_j$ and
$r_{\max}=\max_{1\le j\le q} r_j$. Then
\begin{equation}\label{conv-annulus}
 \sum_{j=1}^q r_jS^1
 =\bigl\{z\in\R^2:\max(0,2r_{\max}-R)\le |z|\le R\bigr\}.
\end{equation}
In particular, this Minkowski sum is convex if and only if
$2r_{\max}\le R$.
\end{lemma}
\begin{proof}
Since the Minkowski sum is compact, connected, and invariant under
rotations about the origin, the set of moduli of its elements is a
closed interval. For $q=1$, this interval is $[r_1,r_1]$. Proceeding by
induction, suppose that the sum of the first $q$ circles is the annulus
with radii $[a,R]$, and add a circle of radius $r>0$. Adding a vector of
length $r$ to a vector of length $t\in[a,R]$ yields all moduli in
$[|t-r|,t+r]$. Thus the maximum modulus is $R+r$ and the minimum modulus
is
\[
 \min_{a\le t\le R}|t-r|
   =\max(0,a-r,r-R).
\]
Substituting $a=\max(0,2r_{\max}-R)$, we get
$\max\{0,2\max(r_{\max},r)-(R+r)\}$, which completes the induction for
\eqref{conv-annulus}. Finally, the set in \eqref{conv-annulus} is convex
if and only if its inner radius $\max(0,2r_{\max}-R)$ is zero.
\end{proof}

\begin{corollary}\label{conv-exceptional-filling}
With the notation of Corollary~\ref{conv-semisimple}, let
\[
 J=\{j:\gfr_j\cong\mathfrak{sl}_2(\R),\ C_j\ne0,
                 \ \operatorname{rank}(L|_{\pfr_j})=2\}.
\]
If $J$ is empty, then $L(K\cdot C)$ is convex. Otherwise, suppose that
there exist an invertible real linear map $T:\R^2\to\R^2$ and positive
numbers $r_j$ ($j\in J$) such that
\[
 L_j(K_j\cdot C_j)=r_jT(S^1)\quad(j\in J),
 \qquad L_j=L|_{\pfr_j}.
\]
If
\begin{equation}\label{conv-polygon}
 2\max_{j\in J}r_j\ \le\ \sum_{j\in J}r_j,
\end{equation}
then $L(K\cdot C)$ is convex.
\end{corollary}
\begin{proof}
As in the proof of Corollary~\ref{conv-semisimple}, we have the
Minkowski sum decomposition $L(K\cdot C)=\sum_j L_j(K_j\cdot C_j)$. By
Theorem~\ref{conv-main}, every summand with $j\notin J$ is convex. For
the summands with $j\in J$, the linearity of $T$,
Lemma~\ref{conv-circle-sums}, and \eqref{conv-polygon} give
\[
 \sum_{j\in J}L_j(K_j\cdot C_j)
   =T\left(\sum_{j\in J}r_jS^1\right)
   =T\left(\{z\in\R^2:|z|\le\sum_{j\in J}r_j\}\right),
\]
which is a closed elliptical disk and hence convex. Therefore,
$L(K\cdot C)$ is a Minkowski sum of convex sets and is convex.
\end{proof}

The assumption that the ellipses are homothetic means that their
principal directions and axis ratios agree. It is stronger than merely
requiring that each projected orbit be an ellipse. Condition
\eqref{conv-polygon} is the familiar polygon inequality: no one radius
exceeds the sum of all the others. The conclusion is a sufficient
condition for a specified pair $(C,L)$, not a necessary condition. In
particular, additional convex summands can fill a hole even when
\eqref{conv-polygon} fails.

\begin{example}\label{conv-so22}
Consider the split real Lie algebra
\[
 \mathfrak{so}(2,2)\cong
 \mathfrak{sl}_2(\R)\oplus\mathfrak{sl}_2(\R).
\]
Let $C=(C_0,C_0)$ with $C_0=\operatorname{diag}(1,-1)$ as in
Remark~\ref{conv-sl2-explanation}, and let
$\phi\bigl(\begin{smallmatrix}x&y\\y&-x\end{smallmatrix}\bigr)=(x,y)$.
For $r_1,r_2>0$, define $L:\pfr_1\oplus\pfr_2\to\R^2$ by
$L(X_1,X_2)=r_1\phi(X_1)+r_2\phi(X_2)$. Then
\[
 L(K\cdot C)=r_1S^1+r_2S^1
 =\{z\in\R^2:|r_1-r_2|\le|z|\le r_1+r_2\}.
\]
Thus $L(K\cdot C)$ is a closed disk (and hence convex) if and only if
$r_1=r_2$, which matches \eqref{conv-polygon}.

By contrast, for
$\mathfrak{so}^*(4)\cong\mathfrak{su}(2)\oplus\mathfrak{sl}_2(\R)$,
the compact simple ideal $\mathfrak{su}(2)$ has trivial $\pfr$-component,
so any rank-two projection of a nontrivial orbit is a single ellipse and is never convex. On the other hand,
$\mathfrak{so}(3,1)\cong\mathfrak{sl}_2(\C)$ is a simple real Lie
algebra not isomorphic to $\mathfrak{sl}_2(\R)$, so all of its planar
orbit projections are convex by Theorem~\ref{conv-main} (see
Helgason~\cite[Chapter~X, Section~6, pages~518--520]{Helgason1978} for
these low-dimensional isomorphisms).
\end{example}

\begin{remark}\label{conv-higher-dim}
The target dimension two in Theorem~\ref{conv-main} is sharp: the
universal convexity of $L(K\cdot C)$ does not extend to real linear maps
$L:\pfr\to\R^n$ for $n\ge3$, even when $\dim\pfr$ is arbitrarily large.
Indeed, consider $\gfr=\mathfrak{sl}_k(\R)$ ($k\ge2$), where
$K=\mathrm{SO}(k)$ acts by orthogonal conjugation on
$\pfr=\{X\in\Sym_k(\R):\tr X=0\}$ of dimension
$\dim\pfr=\frac{k(k+1)}{2}-1$. Take
$C=\diag(k-1,-1,\ldots,-1)\in\pfr$. Every matrix in $K\cdot C$ has the
form $k\,vv^{\top}-I_k$ for a unit vector
$v=(x_1,\ldots,x_k)^{\top}\in S^{k-1}\subset\R^k$. For $X\in\pfr$, define
$L:\pfr\to\R^3$ by 
\[L(X)=(\tr(A_1X),\tr(A_2X),\tr(A_3X)),
\] where
\[
 A_1=\diag(1,-1)\oplus 0_{k-2},\quad
 A_2=\begin{pmatrix}0&1\\1&0\end{pmatrix}\oplus 0_{k-2},\quad
 A_3=\diag(1,1)\oplus 0_{k-2}-\frac{2}{k}I_k.
\]
After adding $(0,0,2)$ and dividing by $k$, the image $L(K\cdot C)$
becomes the set of points
\[
 (u_1,u_2,u_3)=(x_1^2-x_2^2,\,2x_1x_2,\,x_1^2+x_2^2)\in\R^3
 \qquad(v\in S^{k-1}).
\]
Since $u_1^2+u_2^2=u_3^2$, this set is the circle
$\{(u_1,u_2,1):u_1^2+u_2^2=1\}$ when $k=2$, and for every $k\ge3$ it is
the hollow cone surface
\[
 \{(u_1,u_2,u_3)\in\R^3:u_1^2+u_2^2=u_3^2,\ 0\le u_3\le1\},
\]
which omits the interior segment $\{(0,0,u_3):0<u_3\le1\}$ and is
therefore not convex for any $k\ge2$ (cf.~Li and
Tam~\cite[Theorem~4.1, page~150]{LiTam2000}). The complex case behaves differently. For the realification of
$\mathfrak{sl}_k(\C)$, we have $\kfr=\mathfrak{su}(k)$ and
$\pfr=i\,\mathfrak{su}(k)$, on which $K=\mathrm{SU}(k)$ acts by unitary
conjugation. Au-Yeung and Tsing~\cite{AuYeungTsing1983} proved that for
$k\ge3$, \emph{every} three-coordinate linear image of a Hermitian
unitary conjugation orbit is convex. The same is therefore true for
these isotropy orbits; for $k\ge3$, this strengthens Westwick's earlier two-coordinate
$C$-numerical-range theorem~\cite{Westwick1975}. For $k=2$,
a nonzero orbit is a two-sphere, and an isomorphism of $\pfr$ onto
$\R^3$ has a nonconvex image.

Four-coordinate convexity, however, can fail even for $k\ge3$. To see
this, take $C=\diag(k-1,-1,\ldots,-1)$ and write
\[
K\cdot C=\{k\,vv^*-I_k:\|v\|=1\},
\]
 with $v\in\C^k$. Project onto the
four real coordinates
\[
 \bigl(|v_1|^2-|v_2|^2,\,
 2\operatorname{Re}(v_1\overline{v_2}),\,
 2\operatorname{Im}(v_1\overline{v_2}),\,
 |v_1|^2+|v_2|^2\bigr).
\]
This is an affine rescaling of a real linear image of $K\cdot C$. If
the coordinates are $(u_1,u_2,u_3,u_4)$, they satisfy
$u_1^2+u_2^2+u_3^2=u_4^2$ with $0\le u_4\le1$ and hence describe a
nonconvex cone surface. Thus the three-coordinate theorem does not
extend without restriction to four coordinates. More generally,
Li and Poon~\cite{LiPoon2000} studied higher-dimensional joint numerical
ranges as linear images of unitary orbits of rank-$k$ orthogonal
projections. They obtained convexity results under additional
structural conditions and gave examples of nonconvexity with four
coordinates.

There are both dimensional and topological reasons why $n=2$ is the
critical dimension. First, if $C\ne0$ and $L:\pfr\to\R^n$ maps
$\operatorname{span}(K\cdot C)$ onto $\R^n$, then a convex image
$L(K\cdot C)$ must have nonempty interior in $\R^n$. By Sard's theorem
and the fact that $K\cdot C$ lies on the sphere of radius $\|C\|$ in
$\pfr$, this requires the necessary dimension condition
\[
 n\le\dim(K\cdot C)\le\dim\pfr-1.
\]
Second, each step of our proof in Sections~\ref{sec:angular}
and~\ref{sec:zero} relies on $n=2$: in Lemma~\ref{conv-boundary}, an
exposed face of $D=\conv F(M)$ in $\R^n$ lies in a supporting hyperplane
of dimension $n-1$, and a compact connected subset of $\R^{n-1}$ is
automatically convex only when $n-1=1$; in Lemma~\ref{conv-angular}, a
point $p\in\operatorname{int}D\setminus F(M)$ in $\R^n$ defines a
direction map $\theta_p:M\to S^{n-1}$ whose homotopy class lies in
$[M,S^{n-1}]$. Its pullback of the fundamental cohomology class lies
in $H^{n-1}(M;\Z)$, but this cohomology class does not in general
determine the homotopy class; and Lemma~\ref{conv-repzero} applies only
to codimension-two zeros ($L:E\to\R^2$). Consequently, for $n\ge3$,
convexity of $L(K\cdot C)\subset\R^n$ holds only under additional
structural assumptions on $L$ or $C$: for example, when
$n=\dim\mathfrak{a}$ is the real rank of $\gfr$ and
$L:\pfr\to\mathfrak{a}$ is the orthogonal projection onto a maximal
abelian subspace $\mathfrak{a}\subset\pfr$, Kostant's convexity
theorem~\cite{Kostant1973} shows that $L(K\cdot C)=\conv(W\cdot C)$ is a
convex polytope for every $C\in\pfr$; or in the Hermitian case with
$n=3$, when the eigenvalue multiplicities of $C$ and $L$ satisfy the
dimension conditions of Poon~\cite{Poon1997}.

Other three-coordinate convexity results show that the distinction
between two and three coordinates depends on the representation.
In particular, Tam~\cite{Tam1999} proved that three-coordinate
projections of adjoint orbits of $\mathrm{SO}(2m+1)$ on
$\mathfrak{so}(2m+1)$ are convex for $m\ge2$, and gave examples of
nonconvexity with four or more coordinates. These are adjoint orbits of
a compact group rather than the isotropy orbits considered in
Theorem~\ref{conv-main}, so there is no conflict with the general
counterexample above.
\end{remark}

We conclude by comparing Theorem~\ref{conv-main} with related results
in the literature. For compact connected Lie groups, Tam~\cite{Tam2002}
proved the planar convexity of adjoint and coadjoint orbits using the
Kirillov--Kostant--Souriau symplectic structure and Atiyah's connectedness
theorem. Markus and Tam~\cite{MarkusTam2011} subsequently proved
connectedness of fibres of certain maps on compact connected Lie
groups, including the maps $g\mapsto\langle\operatorname{Ad}(g)A,C\rangle$.
Their result gives another proof of the planar convexity of compact
adjoint orbits. These results concern compact-group adjoint orbits,
rather than the isotropy orbits of noncompact symmetric spaces studied
here. For Hermitian unitary conjugation orbits, Westwick
\cite{Westwick1975} established two-coordinate convexity, and Au-Yeung
and Tsing~\cite{AuYeungTsing1983} proved three-coordinate convexity
for matrices of order at least three. Their subsequent paper
\cite{AuYeungTsing1984} also treated generalized numerical ranges
over the real, complex, and quaternionic fields. Poon~\cite{Poon1997}
obtained related results involving eigenvalue-multiplicity hypotheses.
For rotation matrices,
when $n\ge3$, the action of
$\mathrm{SO}(n)\times\mathrm{SO}(n)$ on $M_n(\R)$ given by
$(U,V)X=UXV^{\top}$ is the isotropy representation of
$\mathfrak{so}(n,n)$, and its orbit through $I_n$ is $\mathrm{SO}(n)$.
Thus the convexity of two-dimensional linear projections of
$\mathrm{SO}(n)$ follows as a special case of Theorem~\ref{conv-main}
as well as from Li and Tam~\cite[Theorem~11.7, page~168]{LiTam2000}.
Lau, Ng, and Tsing~\cite{LauNgTsing2016} gave an alternative proof of
planar convexity for special orthogonal orbits and obtained
star-shapedness results for higher-dimensional linear images of these
orbits. More recently, Ramachandran, Shu, and Wang~\cite{RSW2025}
proved the stronger
result that the intersection of $\mathrm{SO}(n)$ with any affine
hyperplane in $M_n(\R)$ is connected, and used this property to design optimization algorithms over
$\mathrm{SO}(n)$. In particular, their hyperplane-section connectedness
is stronger than planar convexity, which was already established in
this case by Li and Tam.

\section*{AI disclosure}
During the preparation of this work, the authors used ChatGPT and Gemini to assist with language editing, exposition, and brainstorm. The authors are responsible for verifying the mathematical arguments and the final manuscript.

\section*{Funding}
Y. Luan is supported by the Shandong Provincial Natural Science Foundation (Grant No. ZR2025MS36).

\end{document}